\documentclass[twocolumn,11pt]{article}
\usepackage{times}
\newtheorem{definition}{Definition}
\newtheorem{remark}[definition]{Remark}

\newtheorem{theorem}[definition]{Theorem}

\newtheorem{corollary}[definition]{Corollary}
\newcommand{\eop}{\hfill $\sqcap\!\!\!\!\sqcup$} 

\usepackage{amssymb,amsmath}
\usepackage{url}

\newtheorem{problem}[definition]{Problem}

\begin{document}
\global\def\refname{{\normalsize \it References:}}
\baselineskip 12.5pt
%
%
%
\title{\LARGE \bf Necessary Optimality Conditions for
Fractional Action-Like Problems with Intrinsic and Observer Times}

\date{}

\author{\hspace*{-10pt}
\begin{minipage}[t]{2.7in} \normalsize \baselineskip 12.5pt
\centerline{GAST\~{A}O S. F. FREDERICO} \centerline{University of
Cape Verde} \centerline{Department of Science and Technology}
\centerline{Praia, Santiago} \centerline{CAPE VERDE}
\centerline{\url{gfrederico@mat.ua.pt}}
\end{minipage} \kern 0in
\begin{minipage}[t]{2.7in} \normalsize \baselineskip 12.5pt
\centerline{DELFIM F. M. TORRES} \centerline{University of Aveiro}
\centerline{Department of Mathematics} \centerline{3810-193
Aveiro} \centerline{PORTUGAL} \centerline{\url{delfim@ua.pt}}
\end{minipage}
\\ \\ \hspace*{-10pt}
\begin{minipage}[b]{6.9in} \normalsize
\baselineskip 12.5pt {\it Abstract:} We prove higher-order
Euler-Lagrange and DuBois-Reymond stationary conditions to
fractional action-like variational problems. More general
fractional action-like optimal control problems are also
considered.
\\ [4mm] {\it Key--Words:}
calculus of variations, FALVA problems, higher-order
Euler-Lagrange equations, higher-order DuBois-Reymond stationary
condition, multi-time control theory.
\end{minipage}
\vspace{-10pt}}

\maketitle


\section{Introduction}
\vspace{-4pt}

The study of fractional problems of the calculus of variations and
respective Euler-Lagrange type equations is a subject of strong
current research because of its numerous applications:
see \textrm{e.g.} \cite{CD:Agrawal:2002,CD:BalAv:2004,CD:Cresson:2007,CD:El-Na:2005,El-Nabulsi2005a,%
tncdf,MR1966935,Klimek2005,CD:Riewe:1996}.
F.~Riewe \cite{CD:Riewe:1996} obtained a version of
the Euler-Lagrange equations for problems of the calculus of
variations with fractional derivatives, that combines both
conservative and non-conservative cases. In 2002 O.~P.~Agrawal proved a
formulation for variational problems with right and left fractional
derivatives in the Riemann-Liouville sense \cite{CD:Agrawal:2002}.
Then, these Euler-Lagrange equations were used by D.~Baleanu and
T.~Avkar to investigate problems with Lagrangians which are linear
on the velocities \cite{CD:BalAv:2004}. In
\cite{MR1966935} fractional problems of the calculus of
variations with symmetric fractional derivatives are considered and
correspondent Euler-Lagrange equations obtained, using both
Lagrangian and Hamiltonian formalisms. In all the above mentioned
studies, Euler-Lagrange equations depend on left and right
fractional derivatives, even when the problem depend only on one
type of them. In \cite{Klimek2005} problems depending on symmetric
derivatives are considered for which Euler-Lagrange equations
include only the derivatives that appear in the formulation of the problem. In \cite{CD:Cresson:2007,El-Nabulsi1}
Euler-Lagrange type equations for problems of the calculus of variations which depend on the
Riemann-Liouville derivatives of order
$\left(\alpha,\beta\right)$, $\alpha > 0$, $\beta > 0$, are studied.

In \cite{CD:DucaU:2006,CD:UdriDu:2005,CD:UdriTe:2007}, C.~Udriste and his coauthors remark that the standard multi-variable variational calculus has
some limitations which the multi-time control theory successfully
overcomes. For instance, the classical multi-variable variational calculus
cannot be applied directly to create a multi-time maximum principle. In \cite{CD:El-Na:2005,El-Nabulsi2005a} two-time Riemann-Liouville fractional integral functionals, depending on a parameter $\alpha$
but not on fractional-order derivatives of order $\alpha$, are
introduced and respective fractional Euler-Lagrange type equations
obtained. In \cite{CD:Jumarie:2007}, Jumarie uses the variational
calculus of fractional order to derive an Hamilton-Jacobi equation,
and a Lagrangian variational approach to the optimal control of
one-dimensional fractional dynamics with fractional cost function. Here, we extend the Euler-Lagrange equations of  \cite{CD:El-Na:2005,El-Nabulsi2005a} by considering two-time fractional action-like variational problems with higher-order derivatives. A DuBois-Reymond stationary condition is also proved for such problems. Finally, we study more general
two-time optimal control type problems.


\section{Preliminaries}
\vspace{-4pt}

In 2005, El-Nabulsi (\textrm{cf.} \cite{CD:El-Na:2005})
introduced the following Fractional Action-Like VAriational (FALVA) problem.
\begin{problem}\label{pb:FAL}
Find the stationary values of the integral functional
\begin{multline} \label{Pi} I[q(\cdot)] =
\frac{1}{\Gamma(\alpha)}\int_a^t
L\left(\theta,q(\theta),\dot{q}(\theta)\right)(t-\theta)^{\alpha-1}
d\theta \tag{$P_1$}
\end{multline}
under given initial conditions
$q(a)=q_{a}$, where $\dot{q} =
\frac{dq}{d\theta}$, $\Gamma$ is the Euler gamma function,
$0<\alpha\leq 1$, $\theta$ is the intrinsic time, $t$ is the observer time,
$t\neq\theta$, and the Lagrangian $L :[a,b] \times
\mathbb{R}^{n} \times \mathbb{R}^{n} \rightarrow \mathbb{R}$ is a
$C^{2}$-function with respect to all its arguments.
\end{problem}

Along all the work, we denote by $\partial_{i}L$ the
partial derivative of a function $L$ with respect to its $i$-th argument, $i \in \mathbb{N}$.

Next theorem summarizes the main result of \cite{CD:El-Na:2005}.

\begin{theorem} (\textrm{cf.} \cite{CD:El-Na:2005})
\label{Thm:NonDtELeq} If $q(\cdot)$ is a solution of
Problem~\ref{pb:FAL}, that is, $q(\cdot)$ offers
a stationary value to functional \eqref{Pi},
then $q(\cdot)$ satisfies the following \emph{Euler-Lagrange equations}:
\begin{multline}
\label{eq:elif}
\partial_{2} L\left(\theta,q(\theta),\dot{q}(\theta)\right)-\frac{d}{d\theta}
\partial_{3} L\left(\theta,q(\theta),\dot{q}(t)\right)\\
= \frac{1-\alpha}{t-\theta}\partial_{3}
L\left(\theta,q(\theta),\dot{q}(\theta)\right)\, .
\end{multline}
\end{theorem}

In this work we begin by generalizing the Euler-Lagrange equations \eqref{eq:elif} for FALVA problems with higher-order derivatives.


\section{Main Results}
\vspace{-4pt}

In \S\ref{sub:sec:EL:FALVA} and \S\ref{sub:sec:DR:FALVA} we
study FALVA problems with higher-order derivatives.
The results are: Euler-Lagrange equations  (Theorem~\ref{Thm:ELdeordm}) and
a DuBois-Reymond stationary condition (Theorem~\ref{theo:cDRifm}) for such problems. Then, on section \S\ref{sub:sec:OC:FALVA}, we introduce the two-time
optimal control FALVA problem, obtaining more general
stationary conditions (Theorems~\ref{theo:pmpif} and \ref{thm:DR:OC}).


\subsection{Euler-Lagrange equations for higher-order FALVA problems} \label{sub:sec:EL:FALVA}
\vspace{-4pt}

We prove Euler-Lagrange equations to higher-order problems of the calculus of variations with fractional integrals of Riemann-Liouville,
\textrm{i.e.} to FALVA problems with higher-order derivatives.

\begin{problem}\label{pb:FALOrdS} The higher-order FALVA problem consists to find stationary values
of an integral functional
\begin{multline}
\label{Pm} I^{m}[q(\cdot)] = \frac{1}{\Gamma(\alpha)}\int_{a}^t
L\Bigl(\theta,q(\theta),\dot{q}(\theta),\\
\ldots,q^{(m)}(\theta)\Bigr)(t-\theta)^{\alpha-1}
d\theta \, , \tag{$P_m$}
\end{multline}
$m\geq 1$, subject to initial conditions
\begin{equation}
\label{eq:Pm} q^{(i)}(a)=q^{i}_{a}\,,\,\, \, i=0, \ldots, m-1\,,
\end{equation}
where
$q^{(0)}(\theta)=q(\theta)$, $q^{(i)}(\theta)$ is the $i$-th derivative, $i\geq
 1$; $\Gamma$ is the Euler gamma function; $0<\alpha\leq 1$; $\theta$ is the intrinsic time;
 $t$ the observer's time, $t\neq\theta$; and the Lagrangian $L :[a,b] \times
\mathbb{R}^{n\times(m+1)} \rightarrow \mathbb{R}$ is a function of
class $C^{2m}$ with respect to all the arguments.
\end{problem}

\begin{remark}
In the particular case
when $m=1$, functional \eqref{Pm} reduces to \eqref{Pi} and Problem~\ref{pb:FALOrdS} to \ref{pb:FAL}.
\end{remark}

To establish the Euler-Lagrange stationary condition
for Problem~\eqref{Pm}, we follow the standard steps used to derive
the necessary conditions in the calculus of variations.

Let us suppose $q(\cdot)$ a solution to
Problem~\ref{pb:FALOrdS}. The variation $\delta I^{m}[q(\cdot)]$
of the integral functional \eqref{Pm} is given by
\begin{equation}
\label{eq:varia}
\frac{1}{\Gamma(\alpha)}\int_{a}^t\left(\sum_{i=0}^m\partial_{i+2}
L\cdot\delta q^{(i)}\right)\left(t-\theta\right)^{\alpha-1}d\theta \, ,
\end{equation}
where $\delta q^{(i)}\in C^{2m}\left([a,b];\mathbb{R}^n\right)$
represents the variation of $q^{(i)}$, $i = 1,\ldots,m$, and satisfy
\begin{equation}\label{eq:cfam}
    \delta q^{(i)}(a)=0\, .
\end{equation}

Having in account conditions \eqref{eq:cfam},
repeated integration by parts of each integral
containing $\delta q^{(i)}$ in
\eqref{eq:varia} leads to
\begin{multline}
\label{eq:elif1} m=1: \quad \delta I[q(\cdot)]
=\frac{1}{\Gamma(\alpha)}\int_{a}^t\Biggl[\left(\partial_{2}
L-\frac{d}{d\theta}\partial_{3} L\right)\\
-\frac{1-\alpha}{t-\theta}\partial_{3}
L\Biggr]\left(t-\theta\right)^{\alpha-1}\cdot\delta q d\theta \, ;
\end{multline}

\begin{multline}
\label{eq:elif2} m=2: \quad \delta I^{2}[q(\cdot)] =
\frac{1}{\Gamma(\alpha)}\int_{a}^t\Biggl[\Bigl(\partial_{2}
L-\frac{d}{d\theta}\partial_{3} L\\
+\frac{d^2}{d\theta^2}\partial_{4}
L\Bigr)
\left(\frac{1-\alpha}{t-\theta}\left(\partial_{3}
L-2\frac{d}{d\theta}\partial_{4}
L\right)\right.\\
\left.-\frac{(1-\alpha)(2-\alpha)}{(t-\theta)^2}\partial_{4}
L\right)\Biggr]\left(t-\theta\right)^{\alpha-1}\cdot\delta q
d\theta \, ;
\end{multline}
and, in general,
\begin{multline*}
\delta I^{m}[q(\cdot)]=\frac{1}{\Gamma(\alpha)}\int_{a}^t\Biggl[\left(\partial_{2}
L+\sum_{i=1}^{m}(-1)^{i}\frac{d^{i}}{d\theta^{i}}\,\partial_{i+2}
L\right)\\
-\frac{1-\alpha}{t-\theta}\sum_{i=1}^{m}i(-1)^{i-1}\frac{d^{i-1}}{d\theta^{i-1}}\,\partial_{i+2} L\\
 - \sum_{k=2}^{m}\sum_{i=2}^{k}(-1)^{i-1}
 \frac{\Gamma(i-\alpha+1)}{(t-\theta)^{i}\Gamma(1-\alpha)}\,
 {k\choose k-i}\frac{d^{k-i}}{d\theta^{k-i}}\,
 \partial_{k+2}L\Biggr]\\
 \cdot \left(t-\theta\right)^{\alpha-1}\cdot\delta q d\theta \, .
\end{multline*}

The integral functional $I^m[\cdot]$ has, by hypothesis, a stationary value for $q(\cdot)$, so that
$$\delta
I^{m}[q(\cdot)]=0\, .$$
The fundamental lemma of the calculus of variations
asserts that all the coefficients
of $\delta q$ must vanish.

\begin{theorem}\textbf{(higher-order Euler-Lagrange equations)}
 \label{Thm:ELdeordm} If $q(\cdot)$ gives a stationary value to functional \eqref{Pm}, then $q(\cdot)$ satisfy the
\emph{higher-order Euler-Lagrange equations}
\begin{multline}
 \label{eq:ELdeordm}
 \sum_{i=0}^{m}(-1)^{i}\frac{d^{i}}{d\theta^{i}}\partial_{i+2}
L\left(\theta,q(\theta),\dot{q}(\theta),\ldots,q^{(m)}(\theta)\right)\\
=F\left(\theta,q(\theta),\dot{q}(\theta),\ldots,q^{(2m-1)}(\theta)\right)\, ,
\end{multline}
where $m\geq 1$ and
\begin{multline}\label{eq:FNCFA}
F\left(\theta,q(\theta),\dot{q}(\theta),\ldots,q^{(2m-1)}(\theta)\right)\\
=\frac{1-\alpha}{t-\theta}\sum_{i=1}^{m}i(-1)^{i-1}\frac{d^{i-1}}{d\theta^{i-1}}\,\partial_{i+2} L\\
\hspace*{-0.3cm}
+\sum_{k=2}^{m}\sum_{i=2}^{k}(-1)^{i-1}
 \frac{\Gamma(i-\alpha+1)}{(t-\theta)^{i}\Gamma(1-\alpha)}\,
 {k\choose k-i}\frac{d^{k-i}}{d\theta^{k-i}}\,
 \partial_{k+2}L
\end{multline}
with the partial derivatives of the Lagrangian $L$ evaluated at $\left(\theta,q(\theta),\dot{q}(\theta),\ldots,q^{(m)}(\theta)\right)$.
\end{theorem}

\begin{remark}
Function $F$ in \eqref{eq:ELdeordm} may be viewed as an external non-conservative friction force acting on the system. If $\alpha=1$, then $F=0$ and equation \eqref{eq:ELdeordm}
is nothing more than the standard Euler-Lagrange equation
for the classical problem of the calculus of variations with
higher-order derivatives:
\begin{equation*}
\sum_{i=0}^{m}(-1)^{i}\frac{d^{i}}{d\theta^{i}}\partial_{i+2}
L\left(\theta,q(\theta),\dot{q}(\theta),\ldots,q^{(m)}(\theta)\right)
=0\, .
\end{equation*}
\end{remark}

\begin{remark}
If $m=1$, the Euler-Lagrange equations \eqref{eq:ELdeordm}
coincide with the Euler-Lagrange equations \eqref{eq:elif}.
\end{remark}

\begin{remark}
For $m=2$, the Euler-Lagrange equations \eqref{eq:ELdeordm} reduce to
\begin{multline}
\label{eq:ELdeord2} \Bigl(\partial_{2}
L\left(\theta,q,\dot{q},\ddot{q}\right)-\frac{d}{d\theta}\partial_{3}
L\left(\theta,q,\dot{q},\ddot{q}\right)\\
+\frac{d^2}{d\theta^2}\partial_{4}
L\left(\theta,q,\dot{q},\ddot{q}\right)\Bigr)
=F\left(\theta,q,\dot{q},\ddot{q},\dddot{q}\right)
\end{multline}
where
\begin{multline}
\label{eq:ELdeord2:aux}
F\left(\theta,q,\dot{q},\ddot{q},\dddot{q}\right)
=\frac{1-\alpha}{t-\theta}\Bigl(\partial_{3}
L\left(\theta,q,\dot{q},\ddot{q}\right)\\
-2\frac{d}{d\theta}\partial_{4}
L\left(\theta,q,\dot{q},\ddot{q}\right)\Bigr)\\
-\frac{\Gamma(3-\alpha)}{(t-\theta)^2\Gamma(1-\alpha)}\,
 {2\choose 0}\partial_{4}
L\left(\theta,q,\dot{q},\ddot{q}\right)\\
=\frac{1-\alpha}{t-\theta}\left(\partial_{3}
L\left(\theta,q,\dot{q},\ddot{q}\right)-2\frac{d}{d\theta}\partial_{4}
L\left(\theta,q,\dot{q},\ddot{q}\right)\right)\\
-\frac{(1-\alpha)(2-\alpha)}{(t-\theta)^2}\partial_{4}
L\left(\theta,q,\dot{q},\ddot{q}\right) \, .
\end{multline}
\end{remark}

\noindent
{\bf Proof:} \
Theorem~\ref{Thm:ELdeordm} is proved by
induction. For $m=1$ and $m=2$, the Euler-Lagrange equations \eqref{eq:elif} and \eqref{eq:ELdeord2}-\eqref{eq:ELdeord2:aux}
are obtained applying the fundamental lemma of the calculus of variations
respectively to \eqref{eq:elif1} and \eqref{eq:elif2}.
From the induction hypothesis,
\begin{multline}
 \sum_{i=0}^{j}(-1)^{i}\frac{d^{i}}{d\theta^{i}}\partial_{i+2}
L\left(\theta,q(\theta),\dot{q}(\theta),\ldots,q^{(m)}(\theta)\right)\\
=F\left(\theta,q(\theta),\dot{q}(\theta),\ldots,q^{(2j-1)}(\theta)\right)
\, ,\, m=j>2.
\end{multline}
We need to prove that equations
\eqref{eq:ELdeordm}--\eqref{eq:FNCFA} hold for $m=j+1$. For simplicity, let us focus our attention on the variation
of $q^{(j+1)}(\theta)$. From hypotheses (variation of
$q^{(j+1)}(\theta)$ up to order $m=j$),
and having in mind that
$C^{j}_{i}+C^{j}_{i+1}=C^{j+1}_{i+1}$ and $m\Gamma(m)=\Gamma(m+1)$, we obtain equations \eqref{eq:ELdeordm} for $m=j+1$
using integration by parts followed by the application
of the fundamental lemma of the calculus of variations.
\eop

\vspace{2pt}

It is convenient to introduce the following quantity (\textrm{cf.} \cite{CD:JMS:Torres:2004}):
\begin{equation}
\label{eq:eqprin}
\hspace*{-0.3cm}
\psi^{j}=\sum_{i=0}^{m-j}(-1)^{i}\frac{d^{i}}{d\theta^{i}}\partial_{i+j+2}
L\left(\theta,q(\theta),\dot{q}(\theta),\ldots,q^{(m)}(\theta)\right),
\end{equation}
$j=1,\ldots,m$. This notation is useful for our purposes because of the following property:
\begin{equation}
\label{eq:eqprin1} \frac{d}{d\theta}\psi^{j}=\partial_{j+1}
L\left(\theta,q(\theta),\dot{q}(\theta),\ldots,q^{(m)}(\theta)\right)-\psi^{j-1},
\end{equation}
$j=1,\ldots,m$.

\begin{remark}
Equation \eqref{eq:ELdeordm} can be written in the following form:
\begin{multline}
\label{eq:ELdeordm1}
\partial_{2}
L\left(\theta,q(\theta),\dot{q}(\theta),\ldots,q^{(m)}(\theta)\right)
-\frac{d}{d\theta}\psi^{1}\\
=F\left(\theta,q(\theta),\dot{q}(\theta),\ldots,q^{(2m-1)}(\theta)\right)\, .
\end{multline}
\end{remark}


\subsection{DuBois-Reymond condition
for higher-order FALVA problems}\label{sub:sec:DR:FALVA}
\vspace{-4pt}

We now prove a DuBois-Reymond condition for FALVA problems.
The result seems to be new even for $m = 1$ (Corollary~\ref{cor:DR:m1}).

\begin{theorem}\textbf{(higher-order DuBois-Reymond condition)} \label{theo:cDRifm} A necessary condition for
$q(\cdot)$ to be a solution to Problem~\ref{pb:FALOrdS} is given by the following higher-order DuBois-Reymond condition:
\begin{multline}
\label{eq:DBRordm}
\hspace*{-0.3cm}
\frac{d}{d\theta}\left\{L\left(\theta,q(\theta),\dot{q}(\theta),\ldots,q^{(m)}(\theta)\right)
-\sum_{j=1}^{m}\psi^{j}\cdot q^{(j)}(\theta)\right\}\\
=\partial_{1}
L\left(\theta,q(\theta),\dot{q}(\theta),\ldots,q^{(m)}(\theta)\right)\\
+F\left(\theta,q(\theta),\dot{q}(\theta),\ldots,q^{(2m-1)}(\theta)\right)\cdot
\dot{q}(\theta)\, ,
\end{multline}
where $F$ and $\psi^j$ are defined by \eqref{eq:FNCFA} and
\eqref{eq:eqprin}, respectively.
\end{theorem}

\begin{remark}
If $\alpha=1$, then $F=0$ and condition \eqref{eq:DBRordm}
is reduced to the classical higher-order DuBois-Reymond condition (see \textrm{e.g.} \cite{CD:JMS:Torres:2004}):
\begin{multline*}
\hspace*{-0.4cm}
\partial_{1} L\left(\theta,q(\theta),\dot{q}(\theta),\ldots,q^{(m)}(\theta)\right)\\
\hspace*{-0.4cm}
=\frac{d}{d\theta}\left\{L\left(\theta,q(\theta),\dot{q}(\theta),\ldots,q^{(m)}(\theta)\right)
-\sum_{j=1}^{m}\psi^{j}\cdot q^{(j)}(\theta)\right\}
\end{multline*}
\end{remark}

\noindent
{\bf Proof:} \
The total derivative of
$$L\left(\theta,q(\theta),\dot{q}(\theta),\ldots,q^{(m)}(\theta)\right)
-\sum_{j=1}^{m}\psi^{j}\cdot q^{(j)}(\theta)$$
with respect to $\theta$ is:
\begin{multline}
\label{eq:FADR}
\frac{d}{d\theta}\left\{L\left(\theta,q(\theta),\dot{q}(\theta),\ldots,q^{(m)}(\theta)\right)
-\sum_{j=1}^{m}\psi^{j}\cdot q^{(j)}(\theta)\right\}\\
=\frac{\partial L}{\partial
\theta}\left(\theta,q(\theta),\dot{q}(\theta),\ldots,q^{(m)}(\theta)\right)\\
+\sum_{j=0}^{m}\frac{\partial L}{\partial
q^{(j)}}\left(\theta,q(\theta),\dot{q}(\theta),\ldots,q^{(m)}(\theta)\right)\cdot
q^{(j+1)}(\theta) \\
-\sum_{j=1}^{m}\left(\dot{\psi}^{j}\cdot
q^{(j)}(\theta)+\psi^{j}\cdot q^{(j+1)}(\theta)\right)\, .
\end{multline}
From \eqref{eq:eqprin1} it follows that \eqref{eq:FADR} is equivalent to
\begin{multline}
\label{eq:FADR1}
\hspace*{-0.3cm}
\frac{d}{d\theta}\left\{L\left(\theta,q(\theta),\dot{q}(\theta),\ldots,q^{(m)}(\theta)\right)
-\sum_{j=1}^{m}\psi^{j}\cdot
q^{(j)}(\theta)\right\}\\
=\frac{\partial L}{\partial
\theta}\left(\theta,q(\theta),\dot{q}(\theta),\ldots,q^{(m)}(\theta)\right)\\
+\sum_{j=0}^{m}\frac{\partial L}{\partial
q^{(j)}}\left(\theta,q(\theta),\dot{q}(\theta),\ldots,q^{(m)}(\theta)\right)\cdot
q^{(j+1)}(\theta)\\
-\sum_{j=1}^{m}\Bigl[\bigl(\frac{\partial L}{\partial
q^{(j-1)}}\left(\theta,q(\theta),\dot{q}(\theta),\ldots,q^{(m)}(\theta)\right)\\
-\psi^{j-1}\bigr)\cdot q^{(j)}(\theta)+\psi^{j}\cdot
q^{(j+1)}(\theta)\Bigr]\, .
\end{multline}
We now simplify the last term on the right-hand side of
\eqref{eq:FADR1}:
 \begin{multline}
\label{eq:FADR2} \sum_{j=1}^{m}\left[\left(\frac{\partial
L}{\partial
q^{(j-1)}}
-\psi^{j-1}\right)\cdot q^{(j)}(\theta)+\psi^{j}\cdot
q^{(j+1)}(\theta)\right]\\
=\sum_{j=0}^{m-1}\Bigl[\frac{\partial L}{\partial
q^{(j)}} \cdot
q^{(j+1)}(\theta) -\psi^{j}\cdot q^{(j+1)}(\theta)\\
+\psi^{j+1}\cdot
q^{(j+2)}(\theta)\Bigr]\\
=\sum_{j=0}^{m-1}\left[\frac{\partial L}{\partial
q^{(j)}}\cdot
q^{(j+1)}(\theta)\right] -\psi^{0}\cdot
\dot{q}(\theta)+\psi^{m}\cdot q^{(m+1)}(\theta)\, ,
\end{multline}
where the partial derivatives of the Lagrangian $L$ are evaluated
at $\left(\theta,q(\theta),\dot{q}(\theta),\ldots,q^{(m)}(\theta)\right)$.
Substituting \eqref{eq:FADR2} into \eqref{eq:FADR1}, and using the higher-order Euler-Lagrange
equations \eqref{eq:ELdeordm}, we obtain the
intended result, that is,
\begin{multline*}
\hspace*{-0.3cm}
\frac{d}{d\theta}\left\{L\left(\theta,q(\theta),\dot{q}(\theta),\ldots,q^{(m)}(\theta)\right)
-\sum_{j=1}^{m}\psi^{j}\cdot q^{(j)}(\theta)\right\}\\
=\frac{\partial L}{\partial
\theta}\left(\theta,q(\theta),\dot{q}(\theta),\ldots,q^{(m)}(\theta)\right)\\
+\frac{\partial L}{\partial
q^{(m)}}\left(\theta,q(\theta),\dot{q}(\theta),\ldots,q^{(m)}(\theta)\right)\cdot
q^{(m+1)}(\theta)\\
+\psi^{0}\cdot \dot{q}(\theta)-\psi^{m}\cdot q^{(m+1)}(\theta)\\
=\partial_{1}
L\left(\theta,q(\theta),\dot{q}(\theta),\ldots,q^{(m)}(\theta)\right)\\
+F\left(\theta,q(\theta),\dot{q}(\theta),\ldots,q^{(2m-1)}(\theta)\right)\cdot
\dot{q}(\theta)\, ,
\end{multline*}
since, by definition,
\begin{equation*}
\psi^{m}=\frac{\partial L}{\partial
q^{(m)}}\left(\theta,q(\theta),\dot{q}(\theta),\ldots,q^{(m)}(\theta)\right)
\end{equation*}
and
\begin{equation*}
\psi^0=
\sum_{i=0}^{m}(-1)^{i}\frac{d^{i}}{d\theta^{i}}\partial_{i+2}
L\left(\theta,q(\theta),\dot{q}(\theta),\ldots,q^{(m)}(\theta)\right)\, .
\end{equation*}
\eop

\vspace{2pt}

\begin{corollary}\textbf{(DuBois-Reymond condition)}
\label{cor:DR:m1}
If $q(\cdot)$ is a solution of
Problem~\ref{pb:FAL}, then the following (first-order) DuBois-Reymond condition holds:
\begin{multline} \label{eq:DBRord1}
\frac{d}{d\theta}\left\{L\left(\theta,q(\theta),\dot{q}(\theta)\right)
-\partial_3 L\left(\theta,q(\theta),\dot{q}(\theta)\right)\cdot \dot{q}(\theta)\right\}\\
=\partial_{1} L\left(\theta,q(\theta),\dot{q}(\theta)\right)
+\frac{1-\alpha}{t-\theta}\partial_3
L\left(\theta,q(\theta),\dot{q}(\theta)\right)\cdot
\dot{q}(\theta)\, .
\end{multline}
\end{corollary}

\noindent
{\bf Proof:} \
For $m=1$, condition \eqref{eq:DBRordm} is reduced to
\begin{multline}
\label{eq:DBRord2}
\frac{d}{d\theta}\left\{L\left(\theta,q(\theta),\dot{q}(\theta)\right)
-\psi^{1}\cdot \dot{q}(\theta)\right\}\\
=\partial_{1}
L\left(\theta,q(\theta),\dot{q}(\theta)\right)
+F\left(\theta,q(\theta),\dot{q}(\theta)\right)\cdot
\dot{q}(\theta)\, .
\end{multline}
Having in mind \eqref{eq:FNCFA} and \eqref{eq:eqprin},
we obtain that
\begin{gather}\label{eq:DBRord3}
\psi^{1}=\partial_3 L\left(\theta,q(\theta),\dot{q}(\theta)\right)\,,\\
\label{eq:DBRord4}
 F\left(\theta,q(\theta),\dot{q}(\theta)\right)
 =\frac{1-\alpha}{t-\theta}\partial_3L\left(\theta,q(\theta),\dot{q}(\theta)\right)\,
.
\end{gather}
One finds the intended equality \eqref{eq:DBRord1}
by substituting the quantities
\eqref{eq:DBRord3} and \eqref{eq:DBRord4} into \eqref{eq:DBRord2}.
\eop

\vspace{2pt}


\subsection{Stationary conditions for optimal control FALVA problems}\label{sub:sec:OC:FALVA}
\vspace{-4pt}

Fractional optimal control problems have been studied in
\cite{Agrawal:2004a,El-Nabulsi1,FDA06}. Here we obtain stationary conditions for two-time FALVA problems
of optimal control. We begin by defining the problem.

\begin{problem}
\label{pb:FALCO}
The two-time optimal control FALVA problem
consists in finding the
stationary values of the integral functional
\begin{equation}
\label{Po1} I[q(\cdot),u(\cdot)] =\frac{1}{\Gamma(\alpha)}
\int_{a}^t
L\left(\theta,q(\theta),u(\theta)\right)(t-\theta)^{\alpha-1}
d\theta \, ,
\end{equation}
when subject to the control system
\begin{equation}
\label{eq:cs:OC:FALVA}
\dot{q}(\theta)=\varphi\left(\theta,q(\theta),u(\theta)\right)
\end{equation}
and the initial condition $q(a)=q_a$.
The Lagrangian $L :[a,b] \times \mathbb{R}^{n}\times
\mathbb{R}^{r} \rightarrow \mathbb{R}$ and the velocity vector
$\varphi:[a,b] \times \mathbb{R}^{n}\times \mathbb{R}^r\rightarrow
\mathbb{R}^{n}$ are assumed to be $C^{1}$ functions with respect to all their
arguments. In accordance with the calculus of variations, we suppose that the control functions $u(\cdot)$
take values on an open set of
$\mathbb{R}^r$.
\end{problem}

\begin{remark}
Problem~\ref{pb:FAL} is a particular case
of Problem~\ref{pb:FALCO} where $\varphi(\theta,q,u)=u$.
FALVA problems of the calculus of variations with higher-order derivatives are also easily written in the optimal control form \eqref{Po1}-\eqref{eq:cs:OC:FALVA}. For example, the integral functional of the second-order FALVA problem of the calculus of variations,
\begin{equation*}
I^2[q(\cdot)] = \frac{1}{\Gamma(\alpha)}\int_{a}^t
L\left(\theta,q(\theta),\dot{q}(\theta),\ddot{q}(\theta)\right)(t-\theta)^{\alpha-1}
d\theta \, ,
\end{equation*}
is equivalent to
\begin{gather*}
\frac{1}{\Gamma(\alpha)}\int_{a}^t
L\left(\theta,q^0(\theta),q^1(\theta),u(\theta)\right)(t-\theta)^{\alpha-1}d\theta
  \, , \\
\begin{cases}
\dot{q}^0(\theta) = q^1(\theta) \, , \\
\dot{q}^1(\theta) = u(\theta) \, .
\end{cases}
\end{gather*}
\end{remark}

We now adopt the Hamiltonian formalism.
We reduce \eqref{Po1}-\eqref{eq:cs:OC:FALVA} to the form \eqref{Pi} by considering the augmented functional:
\begin{multline}
\label{eq:pcond1} J[q(\cdot),u(\cdot),p(\cdot)] \\
= \frac{1}{\Gamma(\alpha)}\int_{a}^t \left[{\cal
H}\left(\theta,q(\theta),u(\theta),p(\theta)\right)
-p(\theta)\cdot\dot{q}(\theta)\right]
d\theta  \, ,
\end{multline}
where the Hamiltonian ${\cal H}$ is defined by
\begin{equation}
\label{eq:Hif} {\cal H}\left(\theta,q,u,p\right)
=L\left(\theta,q,u\right)(t-\theta)^{\alpha-1}+p \cdot
\varphi\left(\theta,q,u\right) \, .
\end{equation}

\begin{definition}\textbf{(Process)}
A pair $(q(\cdot),u(\cdot))$ that satisfies the control system
$\dot{q}(\theta)=\varphi\left(\theta,q(\theta),u(\theta)\right)$ and the initial condition $q(a) = q_a$ of
Problem~\ref{pb:FALCO} is said to be a \emph{process}.
\end{definition}

Next theorem gives the weak Pontryagin maximum principle
for Problem~\ref{pb:FALCO}.

\begin{theorem}
 \label{theo:pmpif}
If $(q(\cdot),u(\cdot))$ is a stationary process for Problem~\ref{pb:FALCO}, then there exists a vectorial function $p(\cdot) \in C^{1}([a,b];\mathbb{R}^{n})$ such that for all $\theta$ the tuple $(q(\cdot),u(\cdot),p(\cdot))$ satisfy the following conditions:
\begin{itemize}
\item the Hamiltonian system
\begin{equation}
\label{eq:Hamif}
\begin{cases}
\dot{q}(\theta)=\partial_{4} {\cal H}(\theta, q(\theta), u(\theta),p(\theta)) \, , \\
\dot{p}(\theta) =-\partial_{2} {\cal H}(\theta, q(\theta),
u(\theta),p(\theta))\, ;
\end{cases}
\end{equation}
\item the stationary condition
\begin{equation}
\label{eq:CEif}
 \partial_{3} {\cal H}(\theta, q(\theta), u(\theta),p(\theta))=0 \, ;
\end{equation}
\end{itemize}
where ${\cal H}$ is given by \eqref{eq:Hif}.
\end{theorem}

\noindent
{\bf Proof:} \
We begin by remarking that the first equation
in the Hamiltonian system, $\dot{q}=\partial_{4} {\cal H}$, is nothing more than the control system \eqref{eq:cs:OC:FALVA}. We write the augmented functional \eqref{eq:pcond1} in the following form:
\begin{equation}
\label{eq:pcond2}
\frac{1}{\Gamma(\alpha)}\int_{a}^t \left[\frac{{\cal
H}-p(\theta)\cdot\dot{q}(\theta)}
{(t-\theta)^{\alpha-1}}\right](t-\theta)^{\alpha-1} d\theta \, ,
\end{equation}
where ${\cal H}$ is
evaluated at $\left(\theta,q(\theta),u(\theta),p(\theta)\right)$.
Intended conditions are obtained
by applying the stationary condition \eqref{eq:elif} to \eqref{eq:pcond2}:
\begin{gather*}
\begin{cases}
\frac{d}{d\theta}\frac{\partial }{\partial \dot{q}}\left[
\frac{{\cal H}-p\cdot\dot{q}}{(t-\theta)^{\alpha-1}}\right]
=\frac{\partial}{\partial q} \left[\frac{{\cal
H}-p\cdot\dot{q}}{(t-\theta)^{\alpha-1}}\right]-
\frac{1-\alpha}{t-\theta}\frac{\partial }{\partial \dot{q}}
\left[\frac{{\cal H}-p\cdot\dot{q}}{(t-\theta)^{\alpha-1}}\right]\\
\frac{d}{d\theta}\frac{\partial }{\partial \dot{u}} \left[
\frac{{\cal H}-p\cdot\dot{q}}{(t-\theta)^{\alpha-1}}\right]
=\frac{\partial}{\partial u}\left[ \frac{{\cal
H}-p\cdot\dot{q}}{(t-\theta)^{\alpha-1}}\right]-\frac{1-\alpha}{t-\theta}\frac{\partial
}{\partial \dot{u}}\left[ \frac{{\cal
H}-p\cdot\dot{q}}{(t-\theta)^{\alpha-1}}\right]\;\,
\end{cases}\\
\Leftrightarrow\;
\begin{cases}
-\dot{p}=\partial_{2} {\cal H}\\
0=\partial_{3} {\cal H}\\
\end{cases}
\end{gather*}
\eop

\begin{remark}
For FALVA problems of the calculus of variations,
Theorem~\ref{theo:pmpif} takes the form of Theorem~\ref{Thm:ELdeordm}.
\end{remark}

\begin{definition}\textbf{(Pontryagin FALVA extremal)} We call any tuple $(q(\cdot),u(\cdot),p(\cdot))$ satisfying Theorem~\ref{theo:pmpif} a \emph{Pontryagin FALVA extremal}.
\end{definition}

Next theorem generalizes the DuBois-Reymond condition
\eqref{eq:DBRordm} to Problem~\ref{pb:FALCO}.

\begin{theorem}
\label{thm:DR:OC}
The following property holds along the
Pontryagin FALVA extremals:
\begin{equation}
\label{eq:PROCO}
\frac{d{\cal H}}{d\theta}(\theta, q(\theta), u(\theta),p(\theta)) =\partial_1
{\cal H}(\theta, q(\theta), u(\theta),p(\theta))\,.
\end{equation}
\end{theorem}

\noindent
{\bf Proof:} \
Equality \eqref{eq:PROCO} is a simple consequence of
Theorem~\ref{theo:pmpif}.
\eop

\vspace{2pt}

\begin{remark}
In the classical framework, \textrm{i.e.} for $\alpha = 1$,
the Hamiltonian ${\cal H}$ does not depend explicitly on $\theta$ when the Lagrangian $L$ and the velocity vector
$\varphi$ are autonomous. In that case, it follows from \eqref{eq:PROCO} that the Hamiltonian ${\cal H}$ (interpreted as energy in mechanics) is conserved.
In the FALVA setting, \textrm{i.e.} for $\alpha\neq 1$, this is no longer true: equality
\eqref{eq:PROCO} holds but we have no conservation of energy since, by definition (\textrm{cf.} \eqref{eq:Hif}), the Hamiltonian ${\cal H}$
is never autonomous (${\cal H}$ always depend explicitly on $\theta$ for $\alpha\neq 1$, thus $\partial_1
{\cal H} \ne 0$).
\end{remark}


\vspace{10pt} \noindent
{\bf Acknowledgements:} \ This work is part
of the first author's PhD project, partially
supported by the \emph{Portuguese Institute for Development}
(IPAD). The authors are also grateful to the support of the
\emph{Portuguese Foundation for Science and Technology} (FCT)
through the \emph{Centre for Research in Optimization and Control}
(CEOC) of the University of Aveiro, cofinanced by the European
Community Fund FEDER/POCI 2010.



\end{document}